\documentclass[11pt,a4paper,twoside,reqno]{amsart}
\usepackage{epsfig,
curves}
\usepackage{float}
\usepackage{amsfonts}
\usepackage{amssymb}

\usepackage{fontenc}
\usepackage{amsmath}
\usepackage{latexsym}

\newcommand{\fnote}[1]{\footnote{\small sharp1}}

\newcommand{\Z}{{\mathbb Z}}
\newcommand{\R}{{\mathbb R}}

\newcommand{\C}{{\mathbb C}}

\newcommand{\M}{\mathbf{M}}

\newcommand{\ZZ}{\mathcal{Z}}
\newcommand{\T}{{\mathbb T}}

\newcommand{\area}{{\rm area}}

\newcommand{\dias}{{\rm dias}}
\newcommand{\sys}{{\rm sys}}
\newcommand{\scg}{{\rm scg}}
\newcommand{\vol}{{\rm vol}}
\newcommand{\id}{{\rm id}}
\newcommand{\length}{{\rm length}}

\newcommand{\cf}{{\it cf.}}
\newcommand{\ie}{{\it i.e.}}

\numberwithin{equation}{section}

\newtheorem{theorem}{Theorem}[section]
\newtheorem{proposition}[theorem]{Proposition}

\newtheorem{lemma}[theorem]{Lemma}

\theoremstyle{definition}
\newtheorem{definition}[theorem]{Definition}

\newtheorem{remark}[theorem]{Remark}

\long\def\forget#1\forgotten{} %

\title[Local extremality of the Calabi-Croke sphere]{Local extremality of the Calabi-Croke sphere for the length of the shortest closed geodesic}

\author[S.~Sabourau]{St\'ephane Sabourau}
\address{
St\'ephane Sabourau - Universit\'e Fran\c{c}ois Rabelais, Tours.
Laboratoire de Math\'ematiques et  Physique Th\'eorique.
CNRS, UMR 6083.
F\'ed\'eration de Recherche Denis Poisson (FR 2964).
Parc de Grandmont, 37400 Tours, France}
\email{sabourau@lmpt.univ-tours.fr}



\begin{document}

\begin{abstract}
Recently, F. Balacheff~\cite{bal} proved that the Calabi-Croke sphere made of two flat $1$-unit-side equilateral triangles glued along their boundaries is a local extremum for the length of the shortest closed geodesic among the Riemannian spheres with conical singularities of fixed area.
We give an alternative proof of this theorem, which does not make use of the uniformization theorem, and extend the result to Finsler metrics.
\end{abstract}

\maketitle


\forget
Rien n'apparait dans le document !
\forgotten

\section{Introduction}

In \cite{cro88}, C. Croke proved that every Riemannian sphere~$S^2$ satisfies
\begin{equation} \label{eq:cro}
\area(S^2) \geq \frac{1}{(31)^2} \, {\scg(S^2)}^2
\end{equation}
where $\scg(S^2)$ represents the length of the shortest closed geodesic on~$S^2$.
This inequality is not sharp and has been improved in \cite{NR02}, \cite{sab04} and \cite{rot06}.
Though the round sphere is a critical point for the quotient $\area/\scg^2$ (and might be a local minimum), \cf~\cite{bal06}, it does not represent the optimal bound.
Actually, it was suggested by E.~Calabi and C.~Croke that the optimal bound is attained by a flat metric with three conical singularities, \cf~\cite{cro88} and~\cite{CK03}.
This surface, called the Calabi-Croke sphere and denoted by~$(S^2,g_0)$, is obtained by gluing two copies of a flat $1$-unit-side equilateral triangle along their boundaries.
The quotient $\area/\scg^2$ of this degenerate metric is equal to~$\frac{1}{2 \sqrt{3}}$. 

In a recent preprint~\cite{bal}, F.~Balacheff proved that the Calabi-Croke sphere is locally extremal  by using the uniformization theorem. \\

The goal of the present article is twofold: 1) to present an alternative (more geometrical) proof of the local extremality of the Calabi-Croke sphere, which does not rely on the uniformization theorem; 2) to establish a Finsler version of it.

A closed geodesic on a Riemannian sphere with conical singularities is a loop which is locally length-minimizing in the neighborhood of every point of~$S^2$ with a conical angle of at least~$2\pi$.
In this definition, there is no condition at a conical singularity of conical angle less than~$2\pi$ (a nontrivial loop passing through such a conical singularity is never length-minimizing) in order to guarantee the existence of a nontrivial closed geodesic on any Riemannian sphere with conical singularities.

Rather than using this notion of a closed geodesic, we will show that as in~\cite{bal}, the Calabi-Croke sphere is locally extremal for the diastole over the one-cycle space.
This diastole, denoted by $\dias_{\ZZ}(S^2)$, is defined by a minimax process over the families of one-cycles sweeping out the sphere, \cf~Section~\ref{sec:bir} or \cite{BS} for a precise definition.
We have
\begin{equation} \label{eq:ds}
\scg(S^2) \leq \dias_{\ZZ}(S^2)
\end{equation}
with equality for the Calabi-Croke sphere, \cf~Section~\ref{sec:bir}.

Thus, showing that the Calabi-Croke sphere is locally extremal for the diastole over the one-cycle space implies that it is locally extremal for the length of the shortest closed geodesic too. \\

In order to specify the topology on the space of metrics considered in this article, let us describe further the Calabi-Croke sphere.

The points $x_{1}$, $x_{2}$ and~$x_{3}$ corresponding to the three conical singularities of angles~$\frac{2\pi}{3}$ of the Calabi-Croke sphere are called \emph{vertices}.
By the uniformization theorem, the sphere~$S^{2}$ has only one conformal structure up to diffeomorphism.
This conformal structure is given by the natural conformal structure on $\C P^{1} \simeq \C \cup \{\infty\}$.
With this identification, from~\cite{tro86,tro90} (see also~\cite{bal}), the Calabi-Croke metric on~$S^{2}$ can be written as
$$
g_{0} = \lambda \, (|z-1| \cdot |z| \cdot |z+1|)^{-\frac{4}{3}} \, |dz|^{2}
$$
for some $\lambda >0$.
Thus, the metric~$g_{0}$ is only defined (and smooth) on \mbox{$S^{2} \setminus \{x_{1},x_{2},x_{3}\}$}.

It is then natural to set the following.
A Riemannian metric with conical singularities~$g$ on~$S^{2}$ is said to be close enough from~$g_{0}$ 
if $g$ is close enough from~$g_{0}$ as a smooth Riemannian metric on~$S^{2} \setminus \{x_{1},x_{2},x_{3}\}$. \\

With this definition, we can now state F.~Balacheff's result~\cite{bal}.

\begin{theorem} \label{theo:main}
Every Riemannian metric with conical singularities~$g$ on~$S^2$ close enough from the Calabi-Croke sphere~$g_0$ satisfies the sharp inequality
$$
\area(S^2,g) \geq \frac{1}{2\sqrt{3}} \, {\dias_{\ZZ}(S^2,g)}^2
$$
with equality if and only if $g$ is homothetic to~$g_0$.
\end{theorem}

From~\eqref{eq:ds}, this result still holds true if one replaces $\dias_{\ZZ}(S^2,g)$ with~$\scg(S^2,g)$. 
Other generalizations are presented in Section~\ref{sec:gen}. \\

Our proof of Theorem~\ref{theo:main} rests on three ingredients:
\begin{enumerate}
\item the construction of a ramified cover from the torus onto the sphere;
\item C.~Loewner's systolic inequality on the torus;
\item Morse theory on the loop space and the one-cycle space of the sphere through minimax principles.
\end{enumerate}

The ramified cover from the torus onto the sphere (allowing us to use C.~Loewner's inequality) was first introduced by the author~\cite[Section~2.5]{sab01} in this context.
It is also the starting point of~\cite{bal}.
While the proof in~\cite{bal} makes use of the uniformization theorem and continues with an analysis of the conformal factor of the metric, requiring a bound on its derivative, our arguments are more geometrical.
In particular, they do not require the uniformization theorem.
Instead, they rely on some geometric observations based on the minimax principle on the loop space and the one-cycle space. \\

Our arguments also carry over to Finsler metrics.
We refer to Section~\ref{sec:finsler} for the necessary definitions regarding Finsler metrics and the following theorem, including the definition of the Finsler Calabi-Croke sphere.

In order to prove the local extremality of the Finsler Calabi-Croke sphere for the diastole over the one-cycle space, we also establish a Finsler version of C.~Loewner's systolic inequality, \cf~Theorem~\ref{theo:loewner}. \\

All this leads to the following.

\begin{theorem} \label{theo:finsler}
Every Finsler metric~$F$ close enough from the Finsler Calabi-Croke sphere~$F_0$ satisfies the sharp inequality
$$
\area(S^2,F) \geq \frac{2}{3 \pi} \, {\dias_{\ZZ}(S^2,F)}^2
$$
for the Holmes-Thompson area (\ie, the canonical symplectic volume of the unit cotangent bundle divided by~$\pi$), with equality if $F$ is homothetic to $F_0$.
\end{theorem}

As above, this result still holds true if one replaces $\dias_{\ZZ}(S^2,F)$ with~$\scg(S^2,F)$. \\

This article is organized as follows.
First, we define the diastole over the loop space and the one-cycle space, \cf~Section~\ref{sec:bir}.
Then, we introduce a degree three ramified cover from the torus onto the sphere in Section~\ref{sec:cover}.
In the sections \ref{sec:P} to~\ref{sec:sys}, we derive geometric estimates on the lengths of some loops on the sphere leading to a comparison between the diastole over the one-cycle space of the sphere and the systole of the torus.
Using C.~Loewner's systolic inequality, we prove Theorem~\ref{theo:main} in Section~\ref{sec:proof}.
Some generalizations are then presented in Section~\ref{sec:gen}. 
Finsler metrics, along with the Finsler Calabi-Croke sphere, are introduced in Section~\ref{sec:finsler}.
In Section~\ref{sec:loewner}, we give an overview of the different proofs of C.~Loewner's inequality in the Riemannian case and establish a Finsler version of it.
The proof of Theorem~\ref{theo:finsler} is then derived.

\section{Closed geodesics and Birkhoff's minimax principle} \label{sec:bir}

The existence of closed geodesics on a smooth Riemannian sphere~$S^2$ has been established by G.~D.~Birkhoff using Morse theory on the free loop space of~$S^2$.
More specifically, consider the free loop space~$\Lambda S^{2}$ formed of piecewise smooth curves $\gamma:S^{1} \to S^{2}$ parametrized proportionnally to arclength and endowed with the compact-open topology.
The subspace of point curves is denoted by~$\Lambda^{0} S^{2}$.

The isomorphism between $\pi_{1}(\Lambda S^{2},\Lambda^{0} S^{2})$ and~$\pi_{2}(S^{2}) \simeq \Z$ allows us to define the diastole over the loop space of~$S^{2}$ by using a minimax principle
$$
\dias_{\Lambda}(S^{2}) = \inf_{(\gamma_{t})} \, \sup_{0 \leq t \leq 1} \length(\gamma_{t})
$$
where $(\gamma_{t})$ runs over all the one-parameter families of loops inducing a generator of~$\pi_{1}(\Lambda S^{2},\Lambda^{0} S^{2}) \simeq \Z$.
Following G.~D.~Birkhoff~\cite{bir}, there exists a closed geodesic of length~$\dias_{\Lambda}(S^{2})$.
In particular, 
$$
\scg(S^{2}) \leq \dias_{\Lambda}(S^{2}).
$$

The existence of a closed geodesic on~$S^2$ can also be proved by using a minimax principle on a different space, namely the one-cycle space with integer coefficients $\ZZ_1(S^2;\Z)$, or $\ZZ_1(S^2)$ for short, endowed with the flat norm topology.
The use of the one-cycle space allows us to cut and paste loops.
We refer to \cite{alm60}, \cite{sab04} or \cite{BS} for a precise definition.
Simply recall that every one-cycle~$z$ decomposes into a denombrable sum of simple loops and that the mass of~$z$, denoted by~$\M(z)$, is the sum of the lengths of these loops counted with multiplicities.
The Almgren isomorphism~\cite{alm60} between $\pi_1(\ZZ_1(S^2;\Z),\{0\})$ and $H_2(S^2;\Z) \simeq \Z$ allows us to define the diastole over the one-cycle space of~$S^2$ as
$$
\dias_{\ZZ}(S^{2}) = \inf_{(z_{t})} \, \sup_{0 \leq t \leq 1} \M(z_{t})
$$
where $(z_{t})$ runs over all the one-parameter families of one-cycles inducing a generator of~$\pi_{1}(\ZZ_1(S^{2}),\{0\})$.
From a result of J.~Pitts~\cite[p.~468]{pit74} (see also~\cite{CC92}), this minimax principle gives rise to a union of closed geodesics (counted with multiplicities) of total length $\dias_{\ZZ}(S^{2})$.
Hence,
$$
\scg(S^{2}) \leq \dias_{\ZZ}(S^{2}).
$$

Of course,
\begin{equation} \label{eq:ZL}
\dias_{\ZZ}(S^{2}) \leq \dias_{\Lambda}(S^{2}).
\end{equation}
Furthermore, equality holds for nonnegatively curved spheres, \cf~\cite{CC92}. \\

The example of a two-sphere of given area with three spikes arbitrarily long shows that the diastole over the loop space~$\dias_{\Lambda}(S^{2})$ provides no lower bound on the area of the sphere (\cf~\cite[Remark 4.10]{sab04} for further details).
However, the diastole over the one-cycle space does provide a lower bound on the area of the sphere similar to C.~Croke's inequality~\eqref{eq:cro}, \cf~\cite{BS}. \\

These definitions extend to Riemannian metrics with conical singularities and Finsler metrics.

\begin{remark} \label{rem:double}
In this article, we do not need the full generality of the one-cycle space arising from geometric measure theory.
We could instead work with the double loop space~$\Gamma$ introduced in~\cite[Appendix]{CC92}, and also used in~\cite{NR02,rot06}.
\end{remark}

\section{A ramified cover from the torus onto the sphere} \label{sec:cover}

The proof of the main theorem relies on the following construction of a ramified cover from the torus onto the sphere.
This ramified cover was introduced by the author~\cite{sab01} in the study of short closed geodesics on Riemannian spheres and used later in~\cite{bal}. \\

By the theory of coverings, there exist a degree three cover $\pi:\T^{2} \to S^{2}$ ramified over the three vertices $x_{1}$, $x_{2}$ and~$x_{3}$ of~$S^{2}$, and a deck transformation map $\rho:\T^{2} \to \T^{2}$ only fixing the ramification points of~$\pi$ with $\rho^{3}=\id_{\T^{2}}$ and $\pi \circ \rho = \pi$.

The ramified cover~$\pi$ can also be constructed in a more geometrical way as follows.
First, cut the sphere along the two minimizing arcs of~$g_{0}$ joining $x_1$ to $x_2$ and $x_1$ to $x_3$. 
This yields a parallelogram with all sides of equal length.
Then, glue three copies of this parallelogram along the two sides between $x_{3}$ and the two copies of~$x_{1}$ to form a hexagon (see the figure below). 
By identifying the opposed sides of this parallelogram, we obtain an equilateral flat torus~$\T^{2}$. 
The isometric rotation, defined on the hexagon, centered at~$x_3$ and permuting the parallelograms, passes to quotient and induces a map $\rho:\T^2 \to \T^2$.
This map gives rise to a degree three ramified cover $\pi:\T^2 \longrightarrow S^2$. \\

\bigskip

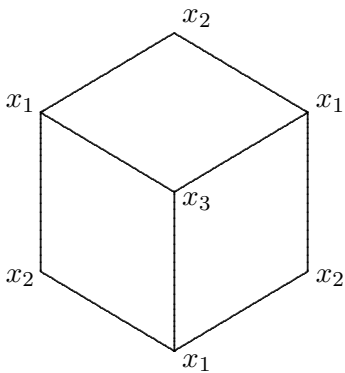
\begin{center}
\setlength\unitlength{1pt}
\begin{picture}(120,140)(-10,-15)

\put(53,-6){$x_1$}
\put(53,124){$x_2$}
\put(-13,26){$x_2$}
\put(103,26){$x_2$}
\put(-13,92){$x_1$}
\put(103,92){$x_1$}
\put(53,55){$x_3$}

\qbezier(0,90)(0,90)(50,120)
\qbezier(0,90)(0,90)(50,60)
\qbezier(0,90)(0,90)(0,30)
\qbezier(100,90)(100,90)(50,120)
\qbezier(100,90)(100,90)(50,60)
\qbezier(100,90)(100,90)(100,30)
\qbezier(50,0)(50,0)(0,30)
\qbezier(50,0)(50,0)(50,60)
\qbezier(50,0)(50,0)(100,30)

\end{picture}
\end{center}

\bigskip 

Thus, the Calabi-Croke sphere can be described as the quotient of an equilateral flat torus by the deck transformation map~$\rho$. \\

Given a Riemannian metric with conical singularities on~$S^{2}$, we will endow~$\T^2$ with the metric pulled back by~$\pi$ and its universal cover~$\R^{2}$ with the metric pulled back by
\begin{equation} \label{eq:cover}
\R^{2} \longrightarrow \T^{2} \stackrel{\pi}{\longrightarrow} S^{2}.
\end{equation}

Since the degree of the Riemannian ramified cover~$\pi$ is equal to three, we have
$$
\area(\T^{2}) = 3 \, \area(S^{2}).
$$

\section{Geometry and dynamics of deformed Calabi-Croke spheres} \label{sec:P}

We will need the following definitions.

\begin{definition} \label{def:sys}
The \emph{systole} of a nonsimply connected metric space~$X$, denoted by~$\sys(X)$, is defined as the infimum of the lengths of the noncontractible loops of~$X$.
The \emph{pointed systole} of~$X$ at a point~$x \in X$, denoted by~$\sys(X,x)$, is defined similarly with the extra assumption that the noncontractible loops considered are based at~$x$.
A shortest noncontractible loop of~$X$ is called a \emph{systolic loop}.
A \emph{systolic loop based at~$x$} is a shortest noncontractible loop of~$X$ with basepoint~$x$.
\end{definition}

As seen in the previous section, the Calabi-Croke sphere can be described as the quotient of an equilateral flat torus by $\pi:\T^{2} \to S^{2}$.
Furthermore, we have the following.

\begin{lemma} \label{lem:close}
A Riemannian metric with conical singularities~$g$ on~$S^{2}$ is close from~$g_{0}$ 
if and only if the pulled-back metric~$\pi^{*}g$ extends to a smooth Riemannian metric on~$\T^2$ close from the flat metric~$\pi^{*}g_{0}$.

In this case, the metric~$g$ has the same conical singularities as $g_{0}$ with the same angles.
\end{lemma}

\begin{proof}
Let $g$ be a metric on~$S^{2}$ given by a function $g:TS^{2} \setminus \Sigma \to \R$ $C^{\infty}$-close from the flat metric~$g_{0}$, with $\Sigma = \cup_{i=1}^3 T_{x_{i}}S^{2}$.
The pulled-back metric~$\pi^{*}g$ on~$\T^{2}$ is also given by a function $\pi^{*}g:\pi^{*}TS^{2} \setminus \pi^{*}\Sigma \to \R$ $C^{\infty}$-close from the flat metric~$\pi^{*}g_{0}$.
Note that $\pi^{*}g_{0}$, contrary to~$\pi^{*}g$, is defined on the whole tangent bundle of~$\T^{2}$.

The metric~$\pi^{*}g$ and all its derivatives are Lipschitz.
From the differential extension theorem, $\pi^{*}g$ extends to a smooth function, still denoted by~$\pi^{*}g$, on~$T\R^{2}$.
The restriction of this function to each tangent plane of~$\R^{2}$ defines a nonnegative quadratic form.
Since $\pi^{*}g$ is close from~$\pi^{*}g_{0}$, this nonnegative quadratic form is nondegenerate.
Therefore, $\pi^{*}g$ defines a smooth metric on~$\R^{2}$, which is close from~$\pi^{*}g_{0}$.
The converse is clear.

For the second part, simply remark that the deck transformation group fixing a given ramification point of~$\pi$ on~$\T^{2}$ is a group isomorphic to~$\Z_{3}$ which acts by isometries on the smooth Riemannian torus~$(\T^{2},\pi^{*}g)$, \cf~\cite{tro90}.
\end{proof}

The Calabi-Croke sphere satisfies the following properties: \\

\begin{enumerate}
\item[(P1)] The length of a shortest loop based at a vertex of~$S^{2}$ which is noncontractible in the sphere with the other two vertices removed is less than twice the minimal distance~$\delta$ between two vertices.
That is, 
$$
\sys(S^{2} \setminus \{x_{i+1},x_{i+2}\},x_{i}) < 2 \delta
$$
for every $i \in \{1,2,3\}$, where the indices are taken modulo~$3$. \\
\item[(P2)] The conical angles of the vertices $x_{1}$, $x_{2}$ and~$x_{3}$ of~$S^2$ are less than~$\pi$. \\
\item[(P3)]  The metric lifts to a smooth Riemannian metric on the universal cover of the torus
$$
\R^{2} \longrightarrow \T^{2} \stackrel{\pi}{\longrightarrow} S^{2},
$$
which does not have any nontrivial closed geodesic of length less than~$3 \, \dias_{\Lambda}(S^{2})$. \\
\end{enumerate}

More generally, we have the following.

\begin{lemma} \label{lem:perturb}
Every Riemannian metric with conical singularities~$g$ on~$S^{2}$ which is close enough from~$g_{0}$ 
still satisfies (P1), (P2) and~(P3).
\end{lemma}

\begin{proof}
The pointed systole and the minimal distance~$\delta$ between two vertices vary continuously with the metric on~$S^{2}$.
Combined with Lemma~\ref{lem:close}, this shows that the open conditions (P1) and~(P2) are satisfied for a metric on~$S^{2}$ close enough from~$g_{0}$.

The metric~$g$ lifts to a smooth Riemannian metric~$\tilde{g}$ on the universal cover of the torus
$$
\R^{2} \longrightarrow \T^{2} \stackrel{\pi}{\longrightarrow} S^{2}.
$$
Observe that the metric~$\tilde{g_{0}}$ defined on~$\R^{2}$ is Euclidean.
The geodesic flow \mbox{$\phi_{\tilde{g}}:T\R^{2} \times \R \to T\R^{2}$} of~$\tilde{g}$ is the solution of a Euler-Lagrange equation given by a second-order differential equation which depends continously on the coefficients of~$\tilde{g}$ and their (first) derivatives, \cf~\cite[Section~1.7]{bes78}.

Let~$\varepsilon>0$.
General results on the dependence of the solutions of a dynamical system with parameters, \cf~\cite{arn06}, imply the following.
If $\tilde{g}$ is close enough from~$\tilde{g_{0}}$, the geodesic flow of~$\tilde{g}$ is $\varepsilon$-close from the geodesic flow of~$\tilde{g_{0}}$ in restriction to any given compact set of $T\R^2 \times \R$.
Furthermore, we can assume that the lengths of the corresponding geodesic trajectories on~$\R^{2}$ are $\varepsilon$-close.
Thus, for $\tilde{g}$ close enough from~$\tilde{g_{0}}$, there is no $\tilde{g}$-geodesic loop of length at most~$3 \, \dias_{\Lambda}(S^2,g)$ on~$\R^{2}$ (remark that the diastole over the loop space -- and the one-cycle space -- of~$S^{2}$ varies continously with the metric on~$S^{2}$).
\end{proof}

\section{curve shortening process} \label{sec:csp}

Consider a Riemannian metric on~$\T^{2}$ which passes to the quotient through $\pi:\T^{2} \to S^{2}$, \cf~Section~\ref{sec:cover}.
The metric on~$\T^{2}$ lifts to a Riemannian metric on the universal cover $\R^{2} \to \T^{2}$.

Let $\gamma$ be a piecewise geodesic loop of~$\R^{2}$ without transverse self-intersection point.
Assume that $\gamma$ is shorter than the shortest closed geodesic on~$\R^{2}$.
The mean curvature flow studied in~\cite{gra89} is a curve shortening process which deforms~$\gamma$ through a one-parameter family of loops~$(\gamma_{t})_{0 \leq t \leq 1}$ where $\gamma_{0}$ agrees with~$\gamma$ and $\gamma_{1}$ is a point curve (after a reparametrization of~$t$).

When $\gamma$ bounds a convex polygonal domain~$D$ (\ie, with acute angles), the mean curvature flow gives rise to a map of degree~$\pm 1$ from the disk onto~$D$.

Furthermore, if $\gamma$ is a lift of a simple loop on~$S^{2}$ through
$$
\R^{2} \longrightarrow \T^{2} \stackrel{\pi}{\longrightarrow} S^{2},
$$
then the same holds for~$\gamma_{t}$ from the equivariance property of the flow.

\begin{remark} \label{rem:csp}
Another equivariant curve shortening process satisfying these properties has been introduced by G.~D.~Birkhoff~\cite{bir}, see also~\cite{cro88}.
To ensure the equivariance property of the flow when $\gamma$ is a lift of a simple loop on~$S^{2}$, one has to take equivariant loop subdivisions in the definition of the curve shortening process.

Furthermore, the Birkhoff curve shortening process extends 
to Finsler metrics, \cf~Section~\ref{sec:finsler}.
\end{remark}

\section{Geodesic loops and diastole} \label{sec:loop}

The following observation is one of the key arguments in the proof of the main theorem.
It compares the length of various geodesic loops on a metric sphere close from the Calabi-Croke sphere.

\begin{lemma} \label{lem:loop}
Let $S^{2}$ be a Riemannian sphere with conical singularities satisfying (P1), (P2) and~(P3), \cf~Section~\ref{sec:P}.

Then, 
$$
\sys(S^{2} \setminus \{x_{i+1},x_{i+2}\},x_{i}) \geq \dias_{\Lambda}(S^{2})
$$
for every $i \in \{1,2,3\}$, where the indices are taken modulo~$3$.
\end{lemma}

\begin{proof}
Let $\gamma$ be a shortest loop based at a vertex of~$S^{2}$ which is noncontractible in the sphere with the other two vertices removed.
Such a loop exists from~(P1).
Furthermore, it is simple and piecewise geodesic (with a singularity of angle less than~$\pi$ only at its basepoint).

We argue by contradiction and suppose that 
$$
\length(\gamma) < \dias_{\Lambda}(S^{2}).
$$
Since the conical angle of the basepoint of~$\gamma$ is less than~$\pi$ from~(P2), the loop~$\gamma$ decomposes~$S^{2}$ into two convex polygonal domains $D_{-}$ and~$D_{+}$.
Lift $D_{-}$ and~$D_{+}$ to two polygonal domains $\tilde{D}_{-}$ and~$\tilde{D}_{+}$ (two triangles) on the universal cover of the torus, \cf~\eqref{eq:cover},
$$
\R^{2} \longrightarrow \T^{2} \stackrel{\pi}{\longrightarrow} S^{2}.
$$
The domains $\tilde{D}_{-}$ and~$\tilde{D}_{+}$ are convex and the length of their boundaries is less than $3 \, \dias_{\Lambda}(S^{2})$.
Applying an equivariant curve shortening process, \cf~Section~\ref{sec:csp}, to the boundaries $\partial \tilde{D}_{-}$ and~$\partial \tilde{D}_{+}$ of the domains $\tilde{D}_{-}$ and~$\tilde{D}_{+}$ induces two one-parameter families of loops from $\partial \tilde{D}_{-}$ and~$\partial \tilde{D}_{+}$ to point curves in $\tilde{D}_{-}$ and~$\tilde{D}_{+}$.
(Indeed, from~(P3), there is no nontrivial closed geodesic on~$\R^{2}$ shorter or of the same length as $\partial \tilde{D}_{-}$ and~$\partial \tilde{D}_{+}$.)
These two families of loops are the lifts of two families of loops $(\gamma^{-}_{t})$ and~$(\gamma^{+}_{t})$ on~$S^{2}$, which both start at~$\gamma$.
Putting together $(\gamma^{-}_{t})$ and~$(\gamma^{+}_{t})$ gives rise to a one-parameter family of loops~$(\gamma_{t})$, starting and ending at point curves, which induces a generator of $\pi_{1}(\Lambda S^{2},\Lambda^{0} S^{2}) \simeq \Z$.
By construction, we have
$$
\length(\gamma) = \sup_{t} \length(\gamma_{t}).
$$
Hence, $\length(\gamma) \geq \dias_{\Lambda}(S^{2})$, which is absurd.
\end{proof}

\section{Projection of systolic loops} \label{sec:proj}

Let $S^{2}$ be a Riemannian sphere with conical singularities and $\pi:\T^{2} \to S^{2}$ be the Riemannian ramified cover of degree three introduced in Section~\ref{sec:cover}.

\begin{lemma} \label{lem:8}
Every systolic loop~$\gamma$ of~$\T^{2}$ which does not pass through a ramification point of~$\pi$ projects to a figure-eight geodesic of~$S^{2}$.

Furthermore, this figure-eight geodesci decomposes~$S^{2}$ into three domains with exactly one vertex lying in each of them.
\end{lemma}

\begin{proof}
The projection of~$\gamma$ on~$S^{2}$ forms a (finite) graph~$\alpha$ with geodesic edges.
Consider the shortest simple loop~$c_{1}$ lying in the support of~$\alpha$ which separates one vertex of~$S^{2}$, say~$x_{1}$ after renumbering, from the other two.
This loop exists, otherwise the lift~$\gamma$ of~$\alpha$ in~$\T^{2}$ would be contractible.
Consider also the shortest simple loop~$c_{2}$ lying in the closure of~$\alpha \setminus c_{1}$ which separates $x_{2}$ from~$x_{3}$.
This loop exists for the same reason as above.
Switching the roles of $x_{2}$ and~$x_{3}$ if necessary, we can assume that the winding number of~$c_{i}$ around the vertex~$x_{j}$ in~$S^{2} \setminus \{x_{3}\}$ is equal to~$\delta_{i,j}$ for $i,j \in \{1,2\}$, where $\delta_{i,j}=1$  if $i=j$ and $0$ otherwise.

Let $c_{3}$ be the shortest arc of~$\alpha$ connecting $c_{1}$ to~$c_{2}$.
We can construct a loop~$c$ made of $c_{1}$, $c_{2}$ and two copies of~$c_{3}$ with winding numbers $1$ and~$-1$ around $x_{1}$ and~$x_{2}$ in~$S^{2} \setminus \{x_{3}\}$.
Clearly, 
\begin{equation} \label{eq:claim}
\length(c) \leq \length(\alpha) = \length(\gamma).
\end{equation}
The loop~$c$ lifts to a noncontractible loop of~$\T^{2}$.
From~\eqref{eq:claim}, we conclude that $c$ is the projection of a systolic loop of~$\T^{2}$ like~$\alpha$.
Thus, $c$ is a geodesic loop of the same length as~$\alpha$.
This implies that $c$ is a figure-eight geodesic loop (with $c_{3}$ reduced to a point) which agrees with~$\alpha$.
Hence the claim.
\end{proof}

\begin{remark}
For a metric on~$S^{2}$ close enough from~$g_{0}$, which is the case we are interested in, we could have used the dynamical properties of the geodesic flow to conclude, \cf~Section~\ref{sec:P}.
\end{remark}

\section{Systole and diastole over the one-cycle space} \label{sec:sys}

In the following result, we compare the diastole over the one-cycle space of some metric spheres with the systole of the corresponding torus.

\begin{proposition} \label{prop:sys}
Let $S^{2}$ be a Riemannian sphere with conical singularities satisfying (P1), (P2) and (P3), and $\pi:\T^2 \to S^2$ be the Riemannian ramified cover of degree three of Section~\ref{sec:cover}.

Then,
$$
\sys(\T^{2}) \geq \dias_{\ZZ}(S^{2}).
$$
\end{proposition}

\begin{proof}
A systolic loop of $S^{2} \setminus \{x_{i+1},x_{i+2}\}$ based at~$x_i$ lifts to a noncontractible loop on the torus.
Thus, its length is greater or equal to the systole of~$\T^2$.
Along with~(P1), we can write
\begin{equation} \label{eq:delta}
\sys(\T^2) \leq \sys(S^{2} \setminus \{x_{i+1},x_{i+2}\},x_i) < 2 \delta.
\end{equation}
Thus, every systolic loop~$\gamma$ of~$\T^{2}$ passes through at most one ramification point of the cover~$\pi$.

If $\gamma$ passes through exactly one ramification point~$y_i$ of the projection~$\pi$, with $\pi(y_{i})=x_{i}$, then its projection is noncontractible in $S^{2} \setminus \{x_{i+1},x_{i+2}\}$.
Otherwise, this projection would lift to a contractible loop of~$\T^{2}$.
From Lemma~\ref{lem:loop}, we obtain
$$
\length(\gamma) \geq \sys(S^{2} \setminus \{x_{i+1},x_{i+2}\},x_{i}) \geq \dias_{\Lambda}(S^2).
$$

If $\gamma$ passes through no ramification point of the projection~$\pi$, then its projection~$\alpha$ is a figure-eight geodesic of~$S^{2}$ which decomposes the sphere into three convex polygonal domains~$D_{i}$, \cf~Lemma~\ref{lem:8}.
As in the proof of Lemma~\ref{lem:loop}, we construct three families of loops from~$\partial D_i$ to point curves in~$D_i$ by using a curve shortening process.
Putting together these three families of loops gives rise to a one-parameter family of one-cycles~$(z_{t})$, starting and ending at null-currents, which induces a generator of $\pi_{1}(\ZZ_1(S^{2}),\{0\}) \simeq \Z$ with
$$
\length(\gamma) = \length(\alpha) = \sup_{t} \M(z_{t}).
$$
Hence,
$
\length(\gamma) \geq \dias_{\ZZ}(S^{2}).
$

Therefore, in both cases, we derive from~\eqref{eq:ZL} that
$$
\sys(\T^{2}) \geq \dias_{\ZZ}(S^{2}).
$$
\end{proof}

\section{Local extremality of the Calabi-Croke sphere} \label{sec:proof}

Before proving Theorem~\ref{theo:main} using the results and the constructions from the previous sections, let us recall C.~Loewner's systolic inequality, \cf~\cite{katz07} for an account on the subject.

Every Riemannian torus~$\T^{2}$ satisfies
\begin{equation} \label{eq:loew}
\area(\T^{2}) \geq \frac{\sqrt{3}}{2} \, {\sys(\T^{2})}^{2}
\end{equation}
with equality if and only if $\T^{2}$ is homothetic to an equilateral flat torus.

We will go over this inequality and extend it to Finsler metrics in Section~\ref{sec:loewner}.

We can now proceed with the proof of Theorem~\ref{theo:main}. \\

Let $S^{2}$ be a Riemannian sphere with conical singularities close enough from the Calabi-Croke sphere to satisfy (P1), (P2) and (P3), \cf~Lemma~\ref{lem:perturb}.
Consider the Riemannian ramified cover of degree three~$\pi:\T^2 \to S^2$ introduced in Section~\ref{sec:cover}.

Since $\area(\T^{2}) = 3 \, \area(S^{2})$, we derive from Proposition~\ref{prop:sys} and C.~Loewner's inequality~\eqref{eq:loew} that
\begin{equation} \label{eq:i}
\area(S^{2}) \geq \frac{1}{2 \sqrt{3}} \, {\dias_{\ZZ}(S^{2})}^{2}.
\end{equation}
Furthermore, if the equality occurs in this latter inequality, then it also occurs in Loewner's.
This implies that $\T^2$ is an equilateral flat torus and that $S^2$ is a sphere obtained by gluing two copies of an equilateral flat triangle along their boundaries.
That is, $S^2$ is homothetic to the Calabi-Croke sphere.

\begin{remark}
The proof of Theorem~\ref{theo:main} shows that the inequality~\eqref{eq:i} holds for every Riemannian sphere with conical singularities satisfying (P1), (P2) and~(P3) (not necessarily close from the Calabi-Croke sphere).
\end{remark}

\section{Some generalizations} \label{sec:gen}

In this section, we briefly discuss some generalizations of Theorem~\ref{theo:main}.
We refrain from giving too much details as the arguments are similar to those presented so far and no new idea is required. \\

As mentionned before, we can replace the diastole over the one-cycle space with the diastole over the double loop space in Theorem~\ref{theo:main}, \cf~Remark~\ref{rem:double}.
Actually, if there is no short loop on~$S^{2}$ representing a local minimum of the length functional, we can replace the diastole over the one-cycle space with the diastole over the loop space in Theorem~\ref{theo:main}.

More specifically, the sharp inequality
$$
\area(S^{2}) \geq \frac{1}{2\sqrt{3}} \, {\min \{ \scg_{0}(S^{2}), \dias_{\Lambda}(S^{2}) \}}^{2}
$$
holds for every Riemannian sphere with conical singularities satisfying (P1) and~(P2), where $\scg_{0}(S^{2})$ is the length of a shortest nontrivial loop of~$S^{2}$ representing a local minimum of the length functional.

Indeed, the right-hand terms in the inequalities of Lemma~\ref{lem:loop} and Proposition~\ref{prop:sys} can be replaced with 
$$
\min \{ \scg_{0}(S^{2}), \dias_{\Lambda}(S^{2}) \}.
$$
This requires some modifications in the proofs but no new idea: simply recall that every loop of length less than $\scg_0(S^2)$ converges to a point through a curve shortening process, and that a systolic loop of~$\T^2$ which does not pass through a ramification point of~$\pi$ projects to a local minimum of the length functional on~$\Lambda S^2$.

One can also show that a Riemannian sphere with conical singularities close enough from the Calabi-Croke sphere for the Lipschitz distance, \cf~\cite{gro99}, satisfies (P1) and (P2) (the angles of the conical singularities vary continously in this topology).
We conclude that the Calabi-Croke sphere is also a local extremum of the quotient $\area/\scg^2$ for the Lipschitz distance topology.

\section{Finsler metrics and the Finsler Calabi-Croke sphere} \label{sec:finsler}

Let us introduce some definitions.

\begin{definition} \label{def:vol}
A function \mbox{$F:TM \longrightarrow \R$}, defined on the tangent bundle of a manifold~$M$, is a Finsler metric if 
\begin{enumerate}
\item it is smooth outside the zero section; \label{1}
\item its restriction on every fiber~$T_x M$ is a vector-space norm; \label{2}
\item $F^2$ has positive definite second derivatives on the fiber~$T_xM \setminus \{0\}$ for every~$x \in M$. \label{3}
\end{enumerate}

One can define the length of piecewise smooth curves with respect to the Finsler metric~$F$ and so the distance between two points as in the Riemannian case.
The condition~\eqref{3} guarantees that the geodesic flow of~$F$ is well defined, \cf~\cite[Chapter~1.F]{bes78}.

The Liouville measure of the cotangent bundle~$T^*M$ of an $n$-dimensional manifold~$M$ is defined as $\frac{\omega^n}{n!}$ where $\omega$ is the canonical symplectic form on~$T^*M$.
The Holmes-Thompson volume of a Finsler $n$-manifold is the Liouville measure of its unit cotangent bundle divided by the Euclidean volume of the unit $n$-ball in $\R^n$.
Unless stated otherwise, we will use this notion of volume on Finsler manifolds.

All these notions, except the existence of a geodesic flow, extend to degenerate Finsler metrics, that is, functions~$F$ satisfying only the conditions \eqref{1} and~\eqref{2}, but not~\eqref{3}.
\end{definition}

\begin{remark} \label{rem:finsler}
The Hausdorff measure provides another way to define the volume of a Finsler manifold or more generally of a metric space.
The Hausdorff measure of an $n$-dimensional normed space~$E$ agrees with the unique Haar measure~$\mu$ such that $\mu(B) = \epsilon_n$, where $B$ is the unit ball of~$E$ and $\epsilon_n$ is the Euclidean volume of the standard unit ball in~$\R^n$, \cf~\cite{tho96}.

Note that the Holmes-Thompson volume and the Hausdorff measure do not necessarily agree.
However, for Riemannian manifolds, the Holmes-Thompson volume and the Hausdorff measure agree with the usual Riemannian volume.
Furthermore, from \cite{dur98}, the Hausdorff measure of a Finsler manifold is not less than its Holmes-Thompson volume, with equality if and only if the metric is Riemannian.
\end{remark}

\begin{definition} \label{def:FCC}
The Riemannian Calabi-Croke sphere~$g_0$ can be described as the quotient of a Riemannian equilateral flat torus by the deck transformation group of the degree three ramified cover $\pi:\T^2 \to S^2$, \cf~Section~\ref{sec:csp}.
Up to homothety, this Riemannian equilateral flat torus agrees with the quotient of the Euclidean plane by the lattice~$\Lambda$ of~$\R^2$ spanned by $\alpha=(1,0)$ and~$\beta=(\frac{1}{2},\frac{\sqrt{3}}{2})$.

Now, consider the Finsler flat torus obtained as the quotient by~$\Lambda$ of the Minkowski plane~$\R^2$ of unit disk the parallelogram of vertices $\frac{\alpha}{2}$,  $\frac{\beta}{2}$, $-\frac{\alpha}{2}$ and~$-\frac{\beta}{2}$.
By analogy with the Riemannian case, we define the Finsler Calabi-Croke sphere~$F_0$ as the quotient of this Finsler flat torus by the deck transformation group of the degree three ramified cover $\pi$.
The metric~$F_0$ is only defined on~$S^2 \setminus \{x_1,x_2,x_3\}$.

A Finsler metric close from the Finsler Calabi-Croke sphere~$F_0$ is a Finsler metric defined on~$S^2 \setminus \{x_1,x_2,x_3\}$ which is close from~$F_0$ as a smooth Finsler metric on $S^2 \setminus \{x_1,x_2,x_3\}$.
\end{definition}

\section{A Finsler Loewner's inequality} \label{sec:loewner}

Sixty years ago, C.~Loewner proved the first systolic inequality, \cf~\cite{katz07}.
Namely, every Riemannian two-torus~$\T^{2}$ satisfies
\begin{equation} \label{eq:loewnerr}
\area(\T^{2}) \geq \frac{\sqrt{3}}{2} \, {\sys(\T^{2})}^{2}
\end{equation}
with equality if and only if $\T^{2}$ is homothetic to a flat equilateral torus.
The original proof of this inequality rests on the uniformization theorem for surfaces, \cf~\cite{katz07}.

There are two other proofs which do not rely on the uniformization theorem.
The first one is a direct consequence of the Burago-Ivanov-Gromov sharp inequality on the $n$-torus for the stable systole, \cf~\cite[Theorem~4.30$_+$]{gro99}.
Simply notice that the systole and the stable systole agree on a Riemannian two-torus.
We will use this approach in the proof of Theorem~\ref{theo:loewner}.

The second one immediately follows from a result of T.~Ilmanen and D.~Knopf.
In~\cite{IK03}, the authors show that the stable systole of a closed Riemannian manifold with nontrivial first homology group is nondecreasing along the Ricci flow.
The same result holds for the systole on an orientable Riemannian surface since, in this case, the systole agrees with the stable systole.
Now, the Ricci flow of a Riemannian two-torus preserves the area and converges to a flat torus, \cf~\cite{ChK04}.
Letting the Ricci flow converge, we derive that the optimal systolic inequality on~$\T^{2}$ is attained by a flat metric.
From Minkowski's theorem, we conclude that this extremal flat torus is equilateral. \\

Now we establish a sharp Finsler version of C.~Loewner's inequality~\eqref{eq:loewnerr}.

\begin{theorem} \label{theo:loewner}
Let $\T^2$ be a Finsler two-torus, then
\begin{equation} \label{eq:loewner}
\area(\T^2)  \geq \frac{2}{\pi} {\sys(\T^2)}^2
\end{equation}
Equality holds if     
$\T^2$ is homothetic to the quotient of $\R^2$, endowed with a parallelogram norm~$||.||$ (i.e., the unit disk of~$||.||$ is a parallelogram), by a lattice with the unit disk of~$||.||$ as a fundamental domain.

The inequality still holds by replacing the Holmes-Thompson area with the Hausdorff measure.
\end{theorem}

\begin{proof}
Let $\T^n$ be a Finsler $n$-torus. The Jacobi torus $J = H_1(\T^n,\R)/H_1(\T^n,\Z)$ is endowed with the stable norm, \cf~\cite[Chapter~4]{gro99}.
Since the Jacobi mapping $f:\T^n \longrightarrow J$ has nonzero degree, we derive from~\cite[Theorem~4.27]{gro99} that 
$$
\vol(\T^n) \geq \vol(J).
$$
Strictly speaking, Theorem~4.27 in~\cite{gro99} states this inequality for a notion of volume different from Holmes-Thompson's.
However, since the Holmes-Thompson volume decreases under distance-decreasing maps, the proof still applies in our case, \cf~\cite[Section~4.29]{gro99}.

We now assume that~$n=2$. Let $D(R)$ be a disk of radius $R < \frac{1}{2} \sys(J)$ in the flat Finsler torus~$J$.
We have 
$
\area(J) \geq \area(D(R)).
$
By definition of the systole, this disk lifts to a disk $\widetilde{D}(R)$ in the normed space~$\widetilde{J}$, where $\widetilde{J}$ is the universal cover of~$J$.
In dimension two, K.~Mahler's inequality, \cf~\cite[Theorem~2.3.4]{tho96} yields the sharp lower bound 
$$
\area(D(R)) = \area(\widetilde{D}(R)) \geq \frac{8}{\pi} R^2
$$
%
where equality holds if and only if $D(R)$ is a parallelogram.
When $R$ tends to $\frac{1}{2} \sys(J)$, we get $\area(J) \geq \frac{2}{\pi} \sys(J)^2$. 
Furthermore, for~$n=2$, we have $\sys(\T^2) = \sys(J)$ since the systole and the stable systole agree in this case.
Therefore, we obtain 
$$
\area(\T^2) \geq \area(J) \geq \frac{2}{\pi} \, {\sys(\T^2)}^2.
$$

If~$\T^2$ is a quotient of~$\R^2$, endowed with a parallelogram norm~$||.||$, by a lattice with the unit disk of~$||.||$ as a fundamental domain, then $\area(\T^2) = \frac{8}{\pi}$ and $\sys(\T^2) = 2$. Hence the result.
\end{proof}

We can now prove Theorem~\ref{theo:finsler}.

\begin{proof}[Proof of Theorem~\ref{theo:finsler}]
The proof follows Theorem~\ref{theo:main}'s, \cf~Section~\ref{sec:proof}, where one replaces C.~Loewner's inequality~\eqref{eq:loewnerr} with its Finsler version, \cf~Theorem~\ref{theo:loewner}.
Simply remark that the intermediate results of the sections \ref{sec:P}, \ref{sec:loop}, \ref{sec:proj} and~\ref{sec:sys} extend to Finsler metrics.
\end{proof}

\section{A Finsler Pu's inequality} \label{sec:pu}

This section is purely expository.
It aims at completing the presentation of sharp systolic inequalities on Finsler surfaces initiated in the previous section.

An analoguous of C.~Loewner's inequality~\eqref{eq:loewner} has been established by P.~Pu~\cite{pu52} on the projective plane.
More precisely, every Riemannian projective plane~$\R P^2$ satisfies 
\begin{equation} \label{eq:pur}
\area(\R P^2) \geq \frac{2}{\pi} \, {\sys(\R P^2)}^2
\end{equation}
with equality if and only if $\R P^2$ is a round projective plane.
As for the torus, the original proof of this inequality rests on the uniformization theorem.

Another proof of this inequality, which does not rely on the uniformization theorem, has been obtained by S.~Ivanov, who showed that \eqref{eq:pur} follows from the main result of~\cite{iva02}.
Furthermore, this proof applies to the Finsler case and yields the same estimate.

Therefore, the estimate of the following theorem is not new. 
We state it as it fits in our study of short closed geodesics on Finsler surfaces.
We simply give the details of the proof and add an observation on the equality case.

\begin{theorem}[see~\cite{iva02}] \label{theo:pu}
Let $\R P^2$ be a Finsler projective plane, then
\begin{equation} \label{eq:pu}
\area(\R P^2) \geq \frac{2}{\pi} \, {\sys(\R P^2)}^2
\end{equation}
Equality holds if 
the geodesic flow of the Finsler metric is periodic or, equivalently, if 
all the geodesics of the Finsler metric are closed.
This is the case, for instance, with the standard round metric.

The inequality still holds by replacing the Holmes-Thompson area with the Hausdorff measure.
In this case, equality holds if and only if the metric is Riemannian with constant curvature.
\end{theorem}

\begin{proof}
Let $\R P^2$ be a Finsler projective plane and $S^2$ be its Finsler universal cover. 
Let $\gamma$ be a loop on $S^2$ which projects onto a systolic loop of~$\R P^2$. 
We have $\length(\gamma) = 2 \, \sys(\R P^2)$. 
The curve $\gamma$ divides~$S^2$ into two isometric disks. 
Let~$D$ be one of them. 
Since $\gamma$ is the lift of a systole, the distance in $D$ between two points of $\gamma$ is given by the length of the shortest subarc of $\gamma$ joining these two points. 
Equivalently, the distance on the boundary of $D$ agrees with the distance on the boundary of the round hemisphere of circumference $\length(\gamma)$, where the boundaries of the two domains are isometrically identified.
From the Finsler version of S.~Ivanov's theorem, \cf~\cite[Section~3]{iva02}, we have 
$$
\area(D) \geq \frac{\length(\gamma)^2}{2 \pi}.
$$ 
Since $\area(D) = \area(\R P^2)$, we finally obtain the Finsler version of P.~Pu's theorem.

Suppose now that the projective plane~$\R P^2$ is endowed with a 
Finsler metric all of whose geodesics are closed. 
The same holds for its universal cover~$S^2$ endowed with the pull-back metric. 
From~\cite{GG81}, the primitive geodesics of~$S^2$ have the same length, namely~$2 \, \sys(\R P^2)$. 
Strictly speaking, the statement of~\cite{GG81} holds for Riemannian metrics, but the proof carries over to Finsler metrics (see also~\cite[p.~143]{zil82}). 
Now, A.~Weinstein's theorem on $C$-manifolds, \cf~\cite[Theorem~2.21 and Proposition~2.24]{bes78}, shows that the Liouville measure of the unit cotangent bundle of~$S^2$ is equal to $4 \, {\sys(\R P^2)}^2$. 
That is, 
$$
\area(\R P^2) = \frac{1}{2} \area(S^2) = \frac{2}{\pi} {\sys(\R P^2)}^2.
$$
The last part of the theorem follows from~\cite{dur98}, \cf~Remark~\ref{rem:finsler}, and the equality case in Pu's inequality~\eqref{eq:pu} on Riemannian projective planes.
\end{proof}

\begin{remark}
Contrary to the Riemannian case, \cf~\cite{bav86}, there is no known sharp systolic inequality on Finsler Klein bottles.
\end{remark}

\end{document}